\documentclass[runningheads]{llncs}

\usepackage{amsmath}
\usepackage{enumerate}
\usepackage[english]{babel}
\usepackage{amssymb}
\usepackage{amscd}
\usepackage{url}
\usepackage{stmaryrd}
\usepackage{tikz,pgf}
\usepackage{subfig}
\usepackage{cancel}
\usepackage{comment}
\usepackage[hidelinks, draft = false]{hyperref}
\usepackage[textsize=footnotesize]{todonotes}
\usepackage{complexity}
\usepackage{cleveref}
\usetikzlibrary{positioning}
\usetikzlibrary{shapes.geometric}
\usepackage{soul}
\usepackage{enumitem}

\usetikzlibrary{arrows}
\usetikzlibrary{decorations.pathreplacing,calligraphy}
\usepgflibrary{patterns.meta} 
\usepgflibrary[patterns.meta] 
\usetikzlibrary{patterns.meta} 
\usetikzlibrary[patterns.meta] 
\definecolor{myBlue}{HTML}{2D2F92}

\setlist[itemize]{noitemsep, topsep=0.5pt}

\usepackage[T1]{fontenc}
\usepackage{graphicx}
\newtheorem{numclaim}[theorem]{Claim}

\newcommand{\np}{\textsf{NP}}

\newcommand{\ie}{{\rm I}_e}
\newcommand{\conf}{{\rm conf}}

\newenvironment{proofclaim}{\noindent{\em Proof of the claim.}}{\qedclaim}

\renewcommand{\qed}{\hfill $\square$ \smallskip}
\renewenvironment{proof}{\noindent{\em Proof.}}{\qed}
\newcommand{\qedclaim}{\hfill $\diamond$ \smallskip}

\newcommand{\defproblem}[3]{
    \vspace{3mm}
    \noindent\fbox{
        \begin{minipage}{0.97\textwidth}
            #1\newline
            \textbf{Input:} #2\\
            \textbf{Question:} #3
        \end{minipage}
    }
    \vspace{3mm}
}

\tikzstyle{v}=[circle,inner sep=0, minimum size =4 pt, line width = 1pt, draw=black, fill=black, text= white]

\begin{document}
\title{Graph Irregularity via Edge Deletions\thanks{This research was partly funded by Kempe Foundation Grant No. JCSMK24-515 (Sweden) and the International Mobility of Researchers MSCA-F-CZ-III at CTU in Prague, $\text{CZ}.02.01.01/00/22\_010/0008601$ Programme Johannes Amos Comenius.}}
%
%
\author{Julien Bensmail\inst{1} \and
Noémie Catherinot\inst{2} \and
Foivos Fioravantes\inst{3} \and 
Clara Marcille\inst{4} \and 
Nacim Oijid\inst{4,5}}

\authorrunning{J. Bensmail et al.}
%

\institute{Universit\'e C\^ote d'Azur, CNRS, Inria, I3S, France \and
Université Paris-Saclay, ENS Paris-Saclay, DER Informatique, 4 Avenue des Sciences 91190 Gif-sur-Yvette, France \and
Dept. of TCS, FIT, Czech Technical University in Prague, Prague, Czech Republic \and
Univ. Bordeaux, CNRS, Bordeaux INP, LaBRI, UMR 5800, F-33400, Talence, France \and 
Department of Mathematics and Mathematical Statistics, Umeå University, Sweden
}


%
\maketitle  

\begin{abstract}
We pursue the study of edge-irregulators of graphs, which were recently introduced in [Fioravantes et al. Parametrised Distance to Local Irregularity. {\it IPEC}, 2024].
That is, we are interested in the parameter ${\rm I}_e(G)$, which, for a given graph $G$,
denotes the smallest $k \geq 0$ such that $G$ can be made locally irregular (\textit{i.e.}, with no two adjacent vertices having the same degree)
by deleting $k$ edges. We exhibit notable properties of interest of the parameter ${\rm I}_e$,
in general and for particular classes of graphs, together with parameterized algorithms for several natural graph parameters. 

Despite the computational hardness previously exhibited by this problem (\NP-hard, \W[1]-hard w.r.t. feedback vertex number, \W[1]-hard w.r.t. solution size), we present two \FPT~algorithms, the first w.r.t. the solution size plus $\Delta$ and the second w.r.t. the vertex cover number of the input graph. 

Finally, we take important steps towards better understanding the behaviour of this problem in dense graphs. This is crucial when considering some of the parameters whose behaviour is still uncharted in regards to this problem (e.g., neighbourhood diversity, distance to clique). In particular, we identify a sub-family of complete graphs for which we are able to provide the exact value of $\ie(G)$. These investigations 
lead us to propose a conjecture that $\ie(G)$ should always be at most $\frac{1}{3}m+c$, where $m$ is the number of edges of the graph $G$ and $c$ is some constant. This conjecture is verified for various families of graphs, including trees. 
\end{abstract}

\section{Introduction}

In graph theory, the family of \textit{regular graphs} is used to capture the most common notion of regularity, that of all vertices having the same degree. 
A very natural question then, is on possible definitions of graphs that would be sort of antonyms to regular graphs.
This led several groups of researchers to investigate various definitions of irregularity in graphs,
such as \textit{irregular graphs} (graphs in which all degrees are pairwise distinct~\cite{CJLORS88}), 
\textit{highly irregular graphs} (graphs in which no vertex has two neighbours with the same degree~\cite{ACCEGO87}), 
and \textit{locally irregular graphs} (graphs in which no two adjacent vertices have the same degree~\cite{BBPW15}).
As we will only focus on the latter notion throughout this work, we encourage the reader to refer to the literature for more details on the former two.

Note that most definitions of (ir)regular graphs do not hold for all simple graphs:
a given simple graph can be more or less (ir)regular w.r.t. some notion of (ir)regularity, and exhibit different behaviour for different such notions.
This leads to different questions, such as how a given simple graph can be turned into a (ir)regular structure through elementary operations. See, e.g.,~\cite{BS21} for an example of a problem where regularity is introduced in the given graph through vertex deletions. We will focus on introducing irregularity. 
The initiators of this line of research are Chartrand \textit{et al.}~\cite{CJLORS88}, who 
introduced the notion of \textit{irregularity strength}.
In brief, this is a parameter that measures how easy/hard it is for a given graph to be turned into an irregular multigraph by increasing the multiplicity of its edges
(with the additional objective that the maximum number of parallel edges replacing a given edge should be minimised).
This concern was considered as well for the other notions of irregularity, namely in~\cite{BLL21} for highly irregular graphs,
and in~\cite{KLT04} for locally irregular graphs. These investigations resulted in numerous related questions and problems,
such as the so-called 1-2-3 Conjecture (introduced by Karo\'nski, {\L}uczak, and Thomason in~\cite{KLT04}, and solved recently 
by Keusch in~\cite{Keu24}) for the latter notion. 

In all these considerations, a given simple graph is intended to be made irregular by essentially adding edges locally (in a parallel way, upon original edges),
and there is an implicit postulate: the number of edges added this way should be minimised.
This leads us to the question of whether one could consider making simple graphs somewhat irregular
through other elementary operations that would limit even further the eventual number of edges.
In particular, could this be achieved without adding any additional edge,
by just \textit{moving} edges? We understand this operation as removing said edge and placing a new edge between any two vertices we chose (regardless of whether there already exists an edge between the considered vertices). At this point it is also worth mentioning~\cite{BCC21}, where the authors consider a similar operation on the edges but with the goal of introducing regularity to the given graph. 

Transforming a given graph into a locally irregular graph via edge moves is straightforward, since any multistar on at least $3$ vertices is always locally irregular.
Things become more interesting when we integrate concerns borrowed from reconfiguration problems, 
such as minimising the number of edge moves performed.
This latter problem is way harder to comprehend. Indeed, when edge moves can be performed at will, it is important to consider the number of edges in a largest (non-induced) locally irregular subgraph of $G$.
This leads us to the study of the (dual) concepts and notions we will focus on in this paper,
which are about the minimum number of edges one must remove from a graph to make it locally irregular. 

These notions were introduced recently by Fioravantes, Melissinos, and Triommatis~\cite{FMT24}.
Following their terminology, we define an \textit{edge-irregulator} of a graph $G$ as a set $S \subseteq E(G)$ of edges of $G$ such that $G-S$ is locally irregular.
We also set $\ie(G)$ as the smallest $k \geq 0$ such that $G$ admits edge-irregulators of cardinality $k$,
and say that an edge-irregulator $S$ of $G$ is \textit{optimal} if $|S|=\ie(G)$.
The authors of~\cite{FMT24} introduced edge-irregulators as an edge counterpart of \textit{vertex-irregulators} (introduced in~\cite{FMT22}),
in which vertices, instead of edges, are deleted. 
It was proven in~\cite{FMT24} that the problem of finding optimal edge-irregulators is \np-hard even when restricted to planar bipartite graphs of maximum degree~$6$,
and \textsf{W[1]}-hard when parameterized by the size of the solution, the feedback vertex set of the input graph, or its vertex integrity.

The first contribution of this work is the continuation of the study of the parameterized complexity of computing an optimal edge-irregulator of a given graph. Considering the previously mentioned \textsf{W[1]}-hardness results, we propose adding the maximum degree $\Delta$ of the input graph as part of the parameter. This proves to be fruitful, as we show that the problem becomes \FPT~w.r.t. the solution size plus $\Delta$. We then present an \FPT~algorithm parameterized by the vertex cover number of the input graph, which is more efficient than the corresponding algorithm that is implied by the work in~\cite{FMT24} (which is based on a reduction to an ILP with bounded number of variables).

The behaviour of computing $\ie$ is less understood when it comes to structural parameters capturing dense graphs. In particular, it is unknown whether the problem is \FPT~ when parameterized by parameters such as the neighbourhood diversity or the distance to clique of the input graph. We take an important step towards understanding the case of complete graphs. In particular, we prove that $\ie(K_{t_k})=|E(K_{t_k})|-\frac{k(k+1)(k-1)(3k+2)}{24}$, where $t_k$ is any triangular number.

The study of this concept of edge-irregulators leads us to conjecture that there exists a constant $c$ such that any connected graph $G$ has an edge-irregulator of size at most $\frac{|E(G)|}{3} + c$. Note that the connected hypothesis is mandatory--the only edge-irregulator of a graph consisting only in a perfect matching will contain all the edges of the graph. Appart form attacking the case of complete graphs, we also prove this conjecture for various families of simple graphs, including trees, for which we also provide a dedicated polynomial algorithm.


\section{Preliminaries, first results and a conjecture}\label{sec:preliminaries}

We follow standard graph theory notations~\cite{D12}. Let $G=(V,E)$ be a graph. For any vertex $v\in V$, let $N_G(v)=\{u\in V : uv\in E\}$ denote the \emph{neighbourhood} of $v$ in $G$, and let $d_G(v)=|N_G(v)|$ denote the \emph{degree} of $v$ in $G$. We also define $N_G[v]=N_G(v)\cup \{v\}$. 
Whenever the graph $G$ is clear from the context, we will omit the subscript and simply write $N(v)$, $d(v)$ and $N[v]$.

Parameterized complexity is a field in algorithm design that takes into consideration additional measures to determine the time complexity. 
The goal in this paradigm is to design a \emph{Fixed-Parameter Tractable} (\FPT) algorithm, i.e., an algorithm that solves the problem in $f(k)|x|^{O(1)}$ time for any arbitrary computable function $f\colon \mathbb{N}\to\mathbb{N}$, where $|x|$ is the size of the input and $k$ is some bounded \textit{parameter}. We say that a problem \textit{is in \FPT} if it admits an \FPT algorithm. 
A parameterized problem is \emph{slicewise polynomial} if it can be determined in $|x|^{f(k)}$ time for a computable function $f\colon \mathbb{N}\to\mathbb{N}$.
Such a problem then belongs to the class \XP.
A problem is presumably not in \FPT if it is shown to be \W[1]-hard (by a parameterized reduction).
We refer the interested reader to classical monographs~\cite{CyganFKLMPPS15,DowneyF13} for a more comprehensive introduction to this topic.

We are now ready to formally define the problem studied in this work. Let $G=(V,E)$ be a graph. We say that $G$ is \emph{locally irregular} if, for every edge $uv\in E$, it holds that $d_G(u)\neq d_G(v)$. Next, let $S\subseteq E$ and consider the graph $G'$ that results from the deletion of the edges of $S$ from $G$. We will denote this $G'$ as $G-S$. Any set $S\subseteq E$ such that $G-S$ is locally irregular will be called called an \emph{edge-irregulator of $G$}. Moreover, let $\ie(G)$ be the minimum order that any edge-irregulator of $G$ can have. We will say that $S$ is an \emph{optimal} edge-irregulator of $G$ if $S$ is an edge-irregulator of $G$ and $|S|=\ie(G)$. Our goal is to compute an optimal edge-irregulator of a given graph $G$ (optimisation version) or, given an integer $k\geq 1$, decide if $\ie(G)\leq k$ (decision version). 

\defproblem{Irregularity-Deletion}{
    A graph $G$, and a positive integer $k \geq 1$. 
}{
    Do we have $\ie(G) \leq k$?
}

To begin, observe that checking if a graph $G$ is locally irregular can be done in polynomial time as it suffices to go through its edges once. Thus, checking if $\ie(G)\leq k$ is (trivially) in \XP~when parameterized by $k$.

Intuitively, regular graphs should be the graph having the largest value of $\ie$. Asymptotically, if $G$ is a cycle with $m$ edges, then all its adjacent vertices are in conflicts and removing an edge resolves at most $3$ of these conflicts. Thus, $\ie(G)$ cannot be too far from $\frac{1}{3}m$.
Equality cannot always hold due to divisibility matters;
for instance, $\ie(C_4)=2$, which is $\frac{1}{2}m$ (where $m=|E(C_4)|=4$) but also $\left\lfloor \frac{1}{3}m \right\rfloor +1$.
Even worse, $\ie(K_2)=1$, which is $m$ (since $m=|E(K_2)|=1$).
Perhaps a good way to express our presumption is in the following form:

\begin{conjecture}\label{conjecture:1/3}
There is an absolute constant $c \geq 1$ such that,
for every connected graph $G$ with $m$ edges, we have $\ie(G) \leq \frac{1}{3}m + c$.
\end{conjecture}

The additive constant $c$ is mandatory to handle some uncontrolled behaviour of small and sparse graphs. In particular, our result on cycles (see Theorem~\ref{theorem:paths}) shows that $c \ge 2$ for any cycle $C_{3k+2}$, where $k$ is an integer.

We want to insist on the fact that the connectivity requirement in Conjecture~\ref{conjecture:1/3} is important,
as the conjecture would be false otherwise. Consider, e.g., a disjoint union of several copies of $K_2$.

Throughout Section~\ref{sec:structure}, we will prove that this conjecture holds for common classes of graphs.
These points apart, Conjecture~\ref{conjecture:1/3} is also supported by the fact that,
when trying to construct connected graphs for which the parameter $\ie$ is a ``large'' fraction of the number of edges,
the best we could achieve is indeed a third of the number of edges. This also shows that Conjecture~\ref{conjecture:1/3} is tight. 

\begin{theorem}\label{theorem:subdiv-1/3}
Let $G$ be any (connected) graph on $m$ edges obtained from a (connected) graph $H$ by subdividing exactly twice every edge of $H$. Then, $\ie(G) \geq \frac{1}{3}m$.
\end{theorem}

\begin{proof}
Consider any path $uvwx$ of $G$, where $d(v)=d(w)=2$, replacing any edge of $H$. Note that, in any edge-irregulator $S$ of $G$, at least one edge in $\{uv,vw,wx\}$ must lie in $S$.
Since $|E(G)|=3|E(H)|$, we thus have $|S| \geq \frac{1}{3}m$.
\end{proof}

As far as lower bounds are concerned, we observe that upon deleting an edge in a graph, the number of \textit{conflicts}, \textit{i.e.}, of remaining edges $uv$ with $d(u)=d(v)$, can decrease by a limited amount only. Formally, for a graph $G$, we define $\conf(G)$ as the number of conflicts in $G$. We stress here the \textit{static} nature of this definition: it only concerns $G$ in its original state. 

\begin{remark}\label{remark:conflict-lower-bound}
In any graph $G$ with maximum degree $\Delta \geq 1$, deleting any edge decreases the number of conflicts by at most $2\Delta-1$. Thus,
$\ie(G)\geq \left\lceil \frac{\conf(G)}{2\Delta-1} \right\rceil $.
\end{remark}

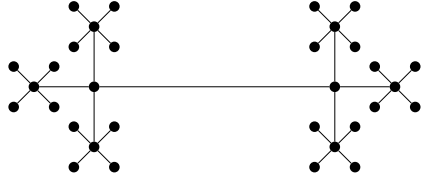
\begin{figure}[t]
    \centering

    \scalebox{0.8}{
    \begin{tikzpicture}
        \draw (0,0) node[v] (c)  {};
        \draw (4,0) node[v] (a)  {};
        \draw (0,1) node[v] (b)  {};
        \draw (0,-1) node[v] (d)  {};
        \draw (-1,0) node[v] (e)  {};
        \draw (5,0) node[v] (f)  {};
        \draw (4,1) node[v] (g)  {};
        \draw (4,-1) node[v] (h)  {};

        \foreach \v in {a, b, d, e} {
            \draw (c) -- (\v);
        }
        \foreach \v in {f,g,h} {
            \draw (a) -- (\v);
        }
        \foreach \v in {b, d, e} {
            \node[v, below right =.3cm of \v] (f1) {};
            \node[v, below left =.3cm of \v] (f2) {};
            \node[v, above right =.3cm of \v] (f3) {};
            \node[v, above left =.3cm of \v] (f4) {};
            
            \draw (\v) -- (f1);
            \draw (\v) -- (f2);
            \draw (\v) -- (f3);
            \draw (\v) -- (f4);
        }
        
        \foreach \v in {f, g, h} {
            \node[v, below right =.3cm of \v] (f1) {};
            \node[v, below left =.3cm of \v] (f2) {};
            \node[v, above right =.3cm of \v] (f3) {};
            \node[v, above left =.3cm of \v] (f4) {};
            
            \draw (\v) -- (f1);
            \draw (\v) -- (f2);
            \draw (\v) -- (f3);
            \draw (\v) -- (f4);
        }
    \end{tikzpicture}
    }
    \caption{Example of graph satisfying $\ie(G) = \left\lceil \frac{\conf(G)}{2\Delta-1} \right\rceil $.}
    \label{fig:conflicts lower bound}
\end{figure}

The bound in Remark~\ref{remark:conflict-lower-bound} is very general, and it can be bad for some graphs.
As an illustration, consider $K_{n,n}$: the bound we obtain is $\ie(K_{n,n})\geq \frac{n^2}{2n-1}$.
Later on, we will show this is distant from the optimal by about a factor of~$2$.
There are other contexts, however, in which this bound is not too bad, as we will see later on with the case of paths and cycles. 
Nevertheless, the bound in Remark~\ref{remark:conflict-lower-bound} can be tight for some graphs, as attested by the graph shown in Figure~\ref{fig:conflicts lower bound}.

\section{FPT~algorithms}\label{sec:algo}

In this section we provide our efficient algorithms for computing an optimal edge-irregulator of a given graph. 

\subsection{Kernelization by solution size plus the maximum degree}

\begin{theorem}\label{thm fpt solution size delta}
    \textsc{Irregularity-Deletion} has a kernel of size $O(k\Delta^{2k+2})$  parameterized by the solution size $k$ plus the maximum degree $\Delta$ of the input graph $G$, and thus an optimal edge-irregulator of $G$ can be computed in \FPT-time.
\end{theorem}

\begin{proof}
    Let $(G,k)$ be an instance of \textsc{Irregularity-Deletion}. Let $\Delta$ be the maximum degree of $G$. For an edge $e \in E(G)$, we denote by $B(e,r)$ the ball of centre $e$ and radius $r$, which is defined as the set of vertices at distance at most $r$ from an endpoint of $e$.

    We first define the set $E_C=\{uv\in E(G):d_G(u)=d_G(v)\}$. That is, $E_C$ contains the edges in conflict of $G$. Recall that, due to Observation~\ref{remark:conflict-lower-bound}, if $|E_C|=\conf(G)\geq k(2\Delta+1)$, then $(G,k)$ is a no-instance of \textsc{Irregularity-Deletion}. Thus, if $|E_C| \ge k(2\Delta+1)$, returns $(K_2, 0)$ as a kernel.

    Otherwise for each $e \in E_C$, we define by $V_e = \{v \in B(e,2k+1) \}$ the set of vertices at distance at most $2k+1$ from an endpoint of $e$. We finally define $V_H = \underset{e \in E_C}{\cup} V_e$ and by $H = G[V_H]$ the subgraph of $G$ induced by $V_H$. Finally, for each vertex $u$ of $V_H$, if $d_H(u) \neq d_G(u)$, we add leaves in $H$ and connect them to $u$ so that it has the same degree in $H$ and in $G$. Denote by $L$ the set of leaves added during this step. Note that this step cannot add conflicts since all the new vertices added have degree one, and since all the vertices of $H$ were already in the connected component of a conflict, they had degree at least one before adding leaves connected to them. Moreover, note by definition, if a vertex is at distance at most $2k-1$ from a conflict $e$, then all its neighbours are in $V_e$. Thus the only vertices on which leaves are added are at distance at least $2k$ from a any conflict. The resulting graph $H$ is a kernel. It remains to prove that $(H,k)$ is a yes instance of  \textsc{Irregularity-Deletion} if and only if $G$ is and to compute the size of $H$.

\begin{numclaim}
    $H$ has at most $ 2 k(2\Delta+1)\Delta^{2k+1}$ vertices
\end{numclaim}
\begin{proofclaim}
    Each $V_e$ is constructed starting from two vertices, and adding at most $\Delta$ new vertices attached to each of them each time the diameter increase. Thus, each $V_e$ contains at most $2 \Delta^{2k}$ vertices. Thus, the union of the $V_e$ contains at most $2 k(2\Delta+1)\Delta^{2k}$ vertices. Finally, we add at most $\Delta$ leaves connected to each of them, which ensure that $|V(H)| \le  2 k(2\Delta+1)\Delta^{2k+1}$.
\end{proofclaim}

\begin{numclaim}
    $(G,k)$ is a yes instance \textsc{Irregularity-Deletion} if and only if $(H,k)$ is.
\end{numclaim}

\begin{proofclaim}
    Suppose first that $(H,k)$ is a yes-instance of \textsc{Irregularity-Deletion}. Let $S$ be an optimal edge-irregulator of $H$. Since all the vertices of $L$ are at distance at least $2k$ from any conflict, using Lemma~\ref{lemma distance edge-irregulator and conflicts}, we know that no edge of $S$ is adjacent to $L$. Thus, all the edges of $S$ are between vertices that are also in $G$. Thus, in $G - S$, if $uv$ is an edge, if both $u$ and $v$ are at distance at most $2k-1$ from a conflict, they have the same degree in $G - S$ and in $H - S$, thus they are not in conflict, and if they are both at distance at least $2k$ from a conflict, no edge of $S$ is adjacent neither to $u$ nor to $v$, thus their degree has not changed, and since they were not in conflict in $G$, they are not in conflict in $G - S$.

    Suppose now that $(G,k)$ is  a yes-instance of \textsc{Irregularity-Deletion}. Let $S$ be an optimal edge-irregulator of $G$. By Lemma~\ref{lemma distance edge-irregulator and conflicts}, all the edges of $S$ are at distance at most $2k-1$ from a conflict, Thus, they are all in $H$. In $H - S$, if $uv$ is an edge, if $u$ or $v$ is in $L$, the other endpoint of $uv$ is at distance exactly $2k+1$ from the nearest conflict. Thus, all its neighbours are at distance at least $2k$ from the nearest conflict, and thus his degree has not changed, thus $uv$ is not a conflict in $H-S$. Otherwise, we have $d_{G-S}(u) = d_{H-S}(u)$ and $d_{G-S}(v) = d_{H-S}(v)$. Thus, since $G - S$ is locally irregular, $uv$ cannot be a conflict.
\end{proofclaim}

Finally, $(G,K)$ and $(H,k)$ are equivalent instances of \textsc{Irregularity-Deletion}, thus $H$ is a kernel and has size $O(k\Delta^{2k+2})$
\end{proof}

The kernel that will be used in Theorem~\ref{thm fpt solution size delta} will be mostly due to the following lemma, and the fact that we can bound the number of vertices in a graph in function of its diameter and its maximum degree. Thus, we prove here the lemma and postpone the proof of Theorem~\ref{thm fpt solution size delta} to the appendix.

\begin{lemma}\label{lemma distance edge-irregulator and conflicts}
    Let $G = (V,E)$ be a graph, and let $S$ be an optimal edge-irregulator of $G$. For any edge $e \in S$, there exists an edge $uv \in E$ such that $d(u) = d(v)$, and $e$ is at distance at most $2|S| - 1$ from $uv$. 
\end{lemma}

\begin{proof}
    The proof is by induction on $k=|S|$.

    Assume first that $k = 1$. Since $S = \{e\}$ is an optimal edge-irregulator and is not empty, there exists at least one conflict $uv \in E$. Thus, since $G - \{e\}$ is locally irregular, either $e = uv$ or one vertex among $u$ or $v$ has a different degree in $G - S$, which implies that $e$ is either adjacent to $u$ or to $v$, and thus at distance at most $1$ from $e$.

    Suppose now $k \ge 2$.
    Let $uv$ be a conflict in $G$. Since in $G-S$, the vertices $u$ and $v$ are not in conflict, there exists one edge $e$ in $S$ that is adjacent to $u$ or to $v$. Consider now the graph $G - \{e\}$. By the optimality of $S$, we get that $S - \{e\}$ is an optimal edge-irregulator of $G - \{e\}$. Thus by induction hypothesis, all the edges of $S$ are at distance at most $2 |S - \{e\}| - 1 = 2 |S| - 3$ from a conflict in $G - \{e\}$. Consider now a conflict $ab$ in $G-\{e\}$. Then either $ab$ was already a conflict in $G$, or $ab$ has become a conflict due to the deletion of $e$. This means that either $a$ or $b$ is adjacent to wither $u$ or to $v$, and thus $ab$ is at distance at most $2$ from $uv$. Finally, using the triangular inequality, any edge of $S$ is at distance at most $2k - 3 + 2 = 2k - 1$ from a conflict of $G$.
\end{proof}

\subsection{Vertex cover number}

While it is already known that $I_e(G)$ can be computed in \FPT-time parameterized by the vertex cover number using the \FPT-algorithm for the more general parameter of \textit{vertex integrity}~\cite{FMT24}, we present here an algorithm that is tailor-made for the vertex cover number, and, as such, is faster. 


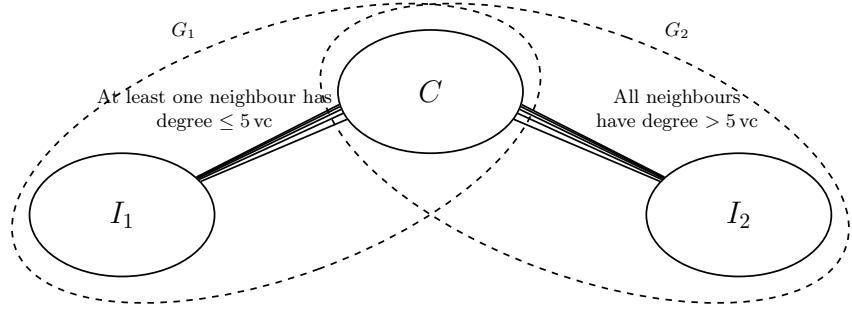
\begin{figure}
    \centering
    \scalebox{0.82}		{
\begin{tikzpicture}[every node/.style={font=\small}]

\tikzset{
    set/.style={draw, thick, ellipse, minimum width=3cm, minimum height=2cm},
    edgestyle/.style={thick},
    labelstyle/.style={font=\bfseries\Large},
    dashedset/.style={draw, ellipse, thick, dashed, inner sep=10pt}
}

\node[set, labelstyle] (I1) at (-2,2) {$I_1$};
\node[set, labelstyle] (I2) at (8,2) {$I_2$};
\node[set, labelstyle] (S)  at (3,4) {$C$};

\foreach \i/\j in {-2/ -2, -1/ -1, 0/ 0, 1/1, 2/2} {
    \draw[edgestyle] ([yshift=1cm,xshift=\i cm]I1) -- ([yshift=-1cm,xshift=\j cm]S);
}

\foreach \i/\j in {-2/ -2, -1/ -1, 0/ 0, 1/1, 2/2} {
    \draw[edgestyle] ([yshift=1cm,xshift=\i cm]I2) -- ([yshift=-1cm,xshift=\j cm]S);
}

\node[align=center] at (-.5,3.7) {At least one neighbour has\\ degree $\leq 5\,\mathrm{vc}$};
\node[align=center] at (7,3.7) {All neighbours\\ have degree $> 5\,\mathrm{vc}$};

\begin{scope}[rotate around={-70:(.5,3)}]
    \draw[dashed, thick] (.5,3) ellipse (2 and 4.5);
\end{scope}
\node at (-1,5) {$G_1$};

\begin{scope}[rotate around={70:(5.5,3)}]
    \draw[dashed, thick] (5.5,3) ellipse (2 and 4.5);
\end{scope}
\node at (7,5) {$G_2$};

\end{tikzpicture}
}
\caption{An abstract depiction of the graph $G$ of Theorem~\ref{thm: FPT vc}.}
    \label{fig:FPT vc}
\end{figure}

\begin{theorem}\label{thm: FPT vc}
    There exists an algorithm to compute $I_e(G)$ in time $2^{O(vc^4)} \cdot  n^{O(1)}$, where $vc$ is the size of a minimum vertex cover of $G$ and $n$ is the order of $G$.
\end{theorem}
\begin{proof}
    We first describe the algorithm and then prove its correctness and time-complexity. 

    Let $G=(V,E)$ be the given graph. Also, let $C=\{u_1,\dots,u_{vc}\}$ be a minimum vertex cover of $G$ and $I=V\setminus C$ be an independent set. We partition $I$ into $I_1$ and $I_2$ as follows. For every vertex $u\in C$, if $d_I(u) \le 5vc$, we include all the vertices of $N_I(u)$ into $I_1$.
    Once we have finished considering all the vertices of $C$, we place the remaining vertices of $I$ into $I_2$. Let $G_1=G[C\cup I_1]$ and $G_2=G[C\cup I_2]$ (see~\Cref{fig:FPT vc}). We will show later that this results in $|I_1|\leq 5vc^2$ and $|V(G_1)|\leq 5vc^2+vc$. We can now guess the projection $S_1$ of an optimal edge-irregulator $S$ of $G$ onto $G_1$. That is, $S_1=S\cap E(G_1)$. Notice that $G_1-S_1$ is not necessarily locally irregular. All that remains to do is to compute the smallest set $S_2\subseteq E(G_2)$ such that $G-(S_1\cup S_2)$ is locally irregular. Notice that such a set $S_2$ contains only edges between the vertices of $C$ and $I_2$; let $E_2$ denote these edges.

    Henceforth, we focus on $G_2$. For every $i\in[vc]=\{1,\dots,vc\}$, we define a \textit{modifier} $k_i\in[n]$. These modifiers encode the number of edges of $E_2$ that need to be removed in order to arrive to a locally irregular graph. That is, if $k_i$ edges of $E_2$ incident to $u_i$ are removed from $G_2$, then $u_i$ is not in conflict with any of its neighbours (in $G$). Note that for each $i$, the modifier $k_i$ does not necessarily have a unique value. We will nevertheless show that it is sufficient to compute a minimum (with respect to summation) tuple $(k_1,\dots,k_{vc})$ such that $\Sigma_{i=1}^{vc}k_i\leq vc^2$. 

    \smallskip 
    
    \noindent\textbf{Correctness.} We first argue about the size of $I_1$.
    \begin{numclaim}\label{claim:sizei1}
        We have that $|I_1|\leq 5vc^2$.
    \end{numclaim}
    \begin{proofclaim}
        Consider how $I_1$ is constructed. By definition of $I_1$, every vertex of $C$ adjacent to a vertex of $I_1$ has degree at most $5vc$. Thus, $|I_1| = |\{N[v] | v \in C, d(v) \le 5vc\}| \le 5vc \cdot |C| = 5vc^2$.
    \end{proofclaim}

    \begin{numclaim}\label{claim:sizee2}
        Let $S_1 \subset E_1$ be a set of edges. For any optimal edge-irregulator $S$ of $G$ coinciding with $S_1$ on $E_1$, we have $S\cap E_2 \leq vc^2$.
    \end{numclaim}
    \begin{proofclaim}
Let $S$ be an optimal edge-irregulator of $G$ coinciding with $S_1$ on $E_1$. Assume by contradiction that $S \cap E_2 > vc ^2$. Since all the edges of $E_2$ has one endpoint in $C$, by the pigeonhole principle, there exists a vertex $v$ in $C$ adjacent to at least $vc +1$ edges of $E_2$. If at least $2vc + 1$ edges of $S$ are adjacent to $v$, then we first remove $vc$ of them from $S$. In this case, $S$ may not yet be an edge-irregulator but the only conflicts we may have created are between $v$ and other vertices of $C$ (since all the other neighbours of $v$ are not in $C$ and thus have degree at most $vc$). Either way, since the other vertices of $C$ can forbid at most $vc - 1$ conflicts, there exist a number $1 \le k \le vc $ such that increasing the degree of $v$ by $k$ removes all these conflicts. Since there are at least $vc+1$ vertices of $S \cap E_2$ that are still adjacent to $v$, we can remove $k$ of them to solve all the conflicts. This leads to a smaller edge-irregulator $S'$ and $S' \cap E_1 = S \cap E_1 = S_1$, which is a contradiction to the optimality of $S$. 
    \end{proofclaim}

  \begin{numclaim}\label{claim:swap}
        Let $S$ be an optimal edge-irregulator of $G$, and $u,v_1,v_2$ be three vertices with $u\in C$ and $v_1,v_2\in I_2$. If $uv_1\in S$, $uv_2\in E$ and $uv_2\notin S$, then there exists an optimal edge-irregulator $S'$ of $G$ such that $S'=(S\setminus \{uv_1\})\cup \{uv_2\}$.
    \end{numclaim}
    \begin{proofclaim}
        First, note that the proof of Claim~\ref{claim:sizee2} also proves that if $S$ is an optimal edge-irregulator of $G$, for any vertex $u \in C$ which has neighbours in $I_2$, we have $d_{G-S}(u)>vc$. 

        Consider now $S'=(S\setminus \{uv_1\})\cup \{uv_2\}$ and let $G'=G-S'$. Since $v_1, v_2 \in I_2$, their degrees are at most $vc$ and all their neighbours have a degree of at least $vc + 1$. Thus $S'$ does not create any conflicts for $v_1$ or for $v_2$. Moreover, all the other vertices in $G$ have the same degree and the same neighbourhood in $G - S$ and in $G - S'$. Thus, since $S$ was an edge-irregulator, they are not in conflict with each other. Thus, $S'$ is an edge-irregulator.
    \end{proofclaim}

    \smallskip 
    
    \noindent\textbf{Running time.} By Claim~\ref{claim:sizei1}, we have that $|C\cup I_1| \leq vc + 5vc^2$. Thus, there are at most $2^{(vc+5vc^2)^2}$ possible guesses of $S_1$.
    By Claim~\ref{claim:sizee2}, we know that an optimal edge-irregulator of $G$ has at most $vc^2$ edges in $E_2$, hence $|S_2|\leq vc^2$. Due to Claim~\ref{claim:swap}, it suffices to guess for each vertex $u\in C$ how many edges of $E_2$ incident to $u$ are in $S_2$, that is, a tuple of integers $(k_1, \dots, k_{vc})$ such that $k_1+\dots +k_{vc}\leq vc^2$. There are $vc^{vc^2}$ such possibilities. Recall that checking whether a set $S\subset E(G)$ is an edge-irregulator can be done in running time $n^2$ where $n$ is the order of $G$, and since we repeat this step for each guess of $S_1$, the running time of this algorithm is $2^{(vc+5vc^2)^2}\cdot vc^{vc^2} \cdot n^2 = 2^{O(vc^4)} \cdot n^{O(1)}$.
\end{proof}

\section{Structural results}\label{sec:structure}



In this section we investigate common classes of graphs to determine their optimal edge-irregulators. Our investigation is lead by Conjecture~\ref{conjecture:1/3}, which we verify for most of the considered classes. The results in this section are important in order to gain a better understanding of $\ie$ in regards to structural parameters, but also since general bounds on $\ie$ can be deduced from simpler cases:

\begin{remark}\label{remark:using-smaller-pieces}
Let $G$ be a graph and $S$ be a set of edges of $G$ such that $G-S$ consists of $p \geq 1$ connected components $G_1,\dots,G_p$. Then, $\ie(G) \leq |S|+\sum_{i=1}^p \ie(G_i)$.
\end{remark}

We start by considering graphs with maximum degree~$2$ (paths and cycles).
For $n \geq 2$, we denote by $P_n$ the path of order~$n$, and, for $n \geq 3$, by $C_n$ the cycle of order~$n$.
In this case, we are able to determine $\ie$ completely. Roughly, to prove this result we delete the minimum number of edges to reach as many disjoint $P_3$s as possible. 

\begin{theorem}\label{theorem:paths}
For every $n \geq 2$, we have:
\begin{itemize}
    \item $\ie(P_n) = \frac{n-1}{3} $ if $n \equiv 1 \bmod 3$;
    \item $\ie(P_n) = \lceil \frac{n-1}{3} \rceil$ if $n \equiv 2 \bmod 3$;
    \item $\ie(P_n) = \lfloor \frac{n-1}{3} \rfloor$ otherwise.
\end{itemize}
Moreover, for every $n \geq 3$, we have $\ie(C_n)=\ie(P_n)+1$.
\end{theorem}

\begin{proof}
Note that the only value of $n \geq 2$ such that $P_n$ is locally irregular is $n=3$. Thus, for any set $S \subseteq E(P_n)$, for $P_n-S$ to be locally irregular, every non-trivial connected component must contain exactly two edges. In particular, $P_3$ is locally irregular, so $\ie(P_3)=0$.

Let us denote by $e_1,e_2,\dots,e_{n-1}$ the consecutive edges of $P_n$, from one end-vertex to the other. Set $E_1=\{e_1,e_2,e_3\}$, $E_2=\{e_4,e_5,e_6\}$, and so on; that is, let $E_1,\dots,E_p$ be consecutive groups of three consecutive edges of $P_n$, where the last group, $E_p$, might contain one (if $n-1 \equiv 1 \bmod 3$) or two (if $n-1 \equiv 2 \bmod 3$) edges only. Note that $p=\lceil \frac{n-1}{3} \rceil$. For lower bounds, observe that, for any set $S \subseteq E(P_n)$ such that $P_n-S$ is locally irregular, at least one edge of each of $E_1,\dots,E_{p-1}$ must belong to $S$. This is because, otherwise, $P_n-S$ would have a connected component containing at least three edges, a contradiction as mentioned earlier. Now, if $n-1 \equiv 0 \bmod 3$, then $p=\frac{n-1}{3}$, and $|S| \geq p$.
If $n-1 \equiv 2 \bmod 3$, then $|S| \geq \lfloor \frac{n-1}{3} \rfloor$. It remains to consider when $n-1 \equiv 1 \bmod 3$, in which case we similarly obtain $|S| \geq \lfloor \frac{n-1}{3} \rfloor$ only. We claim $S$ must contain at least one additional edge. If $S$ contains $e_{n-1}$ (the only edge of $E_p$), then we have our conclusion. Otherwise, by earlier arguments we deduce that $S$ must contain $e_{n-3},e_{n-6}, \dots$, \textit{i.e.}, the second edge of each of $E_1,\dots,E_{p-1}$. If $S$ does not contain any other edge, then we deduce that $e_2$ is part of $S$ while $e_1$ is not, which contradicts that $P_n-S$ is locally irregular. Thus, we must have $|S| \geq \lceil \frac{n-1}{3} \rceil$.

Regarding upper bounds, let us consider $S$, the set of edges of $P_n$ comprising every three edges of $P_n$, namely $e_3,e_6,e_9,\dots$. If $n-1 \equiv 0 \bmod 3$, then note that $e_{n-1}$ belongs to $S$. If $n-1 \equiv 2 \bmod 3$, then the last edge of $P_n$ in $S$ is $e_{n-3}$. When $n-1 \equiv 1 \bmod 3$, the last edge in $S$ is $e_{n-2}$; in that special case, we also add $e_{n-1}$ to $S$. In all three cases, it can be checked that $P_n-S$ is locally irregular, essentially because for every three consecutive edges of $P_n$ there is one in $S$. Now:

\begin{itemize}
    \item if $n-1 \equiv 0 \bmod 3$, then $|S|=\frac{n-1}{3}$;
    \item if $n-1 \equiv 1 \bmod 3$, then $|S|=\lfloor \frac{n-1}{3} \rfloor +1 = \lceil \frac{n-1}{3} \rceil$;
    \item if $n-1 \equiv 2 \bmod 3$, then $|S|= \lfloor \frac{n-1}{3} \rfloor$.
\end{itemize}

\noindent Thus, we have the claimed upper bounds, and all equalities hold.

If $C_n$ is a cycle, all its vertices have degree~$2$, thus, at least one edge has to be removed from $C_n$. After this first edge removed, we are left to a path on $n$ vertices, which concludes the proof.
\end{proof}

Next, we deal with complete bipartite graphs. Recall that $K_{n,m}$ denotes the complete bipartite graph with parts of size $n$ and $m$ ($n,m \geq 1$).

\begin{theorem}\label{theorem:complete-bipartite}
For every $n,m \geq 1$, we have:
\begin{itemize}
    \item $\ie(K_{n,m})=0$ if $n \neq m$;
    \item $\ie(K_{n,m})=n$ otherwise.
\end{itemize}
\end{theorem}

\begin{proof}
The first item is because $K_{n,m}$ is locally irregular whenever $n \neq m$. We now focus on $K_{n,n}$ with $n \geq 2$. Choose any vertex, $u$, of $K_{n,n}$, and let $S$ be the set containing all edges incident to $u$. Then, clearly, $K_{n,n}-S$ is locally irregular, and $\ie(K_{n,n}) \leq n$. To see, now, that $\ie(K_{n,n}) \geq n$, consider, towards a contradiction, any set $S$ of at most $n-1$ edges such that $K_{n,n}-S$ is locally irregular. Denoting by $(U,V)$ the bipartition of $K_{n,n}$ (where, thus, $|U|=|V|=n$), since $|S| \leq n-1$, there must be at least a $u \in U$ not incident to any edge of $S$, and similarly a $v \in V$ not incident to any edge of $S$. Then, $d_{K_{n,n}-S}(u)=d_{K_{n,n}-S}(v)=n$, while $u$ and $v$ are adjacent in $K_{n,n}-S$; this is a contradiction.
\end{proof}

We now move towards more general families of graphs. 

\begin{theorem}\label{theorem:bipartite-deg1-upper}
If $G=(V,E)$ is a connected bipartite graph of order $n$ and minimum degree~$1$, then $\ie(G) \leq n-1$.
\end{theorem}

\begin{proof}
Let $r$ be any degree-$1$ vertex of $G$, and let $V_1,V_2,\dots$ be the partition of $V(G) \setminus \{r\}$ where each $V_i$ contains the vertices of $G$ at distance exactly $i$ from $r$. Our goal is to construct a set $S$ of edges of $G$ such that, in $G-S$, the vertices in some $V_i$ with $i$ even have even degrees, while those in some $V_i$ with $i$ odd have odd degrees (except maybe for vertices becoming isolated, and one vertex; see later). Towards this, we start with $S$ being empty. Assuming $V(G) \setminus \{r\}=V_1 \cup \dots \cup V_d$ (\textit{i.e.}, the maximum distance between $r$ and a vertex of $G$ is $d$), we consider the vertices of $V(G) \setminus \{r\}$ one by one, following the parts $V_d,V_{d-1},\dots,V_1$ in this order. For any vertex $u \neq r$ considered this way, we consider $d_{G-S}(u)$, and, in case this degree does not have the required parity, we add to $S$ a single edge $uu'$ where $u'$ is closer to $r$ than $u$ is ($u'$ exists since $u \neq r$). As a result, after having treated $u$, note that $d_{G-S}(u)$ (for the updated $S$) is as desired modulo~$2$.

Once all vertices of $V_d,\dots,V_1$ have been treated this way, note that, for the resulting $S$, we have that $d_{G-S}(u)$ has the desired parity for every $u \in V(G) \setminus \{r\}$. In particular, for any two adjacent vertices $u$ and $u'$ such that $r \not \in \{u,u'\}$, we have $d_{G-S}(u) \not \equiv d_{G-S}(u') \bmod 2$. So the only reason why $G-S$ might not be locally irregular, is because $r$ has the same (non-zero) degree as its unique neighbour, $w$. In particular, this means $rw \not \in S$, and $r$ is the only neighbour of $w$ in $G-S$. In that case, we simply add $rw$ to $S$ to get rid of this conflict, without creating new ones; then, for the resulting $S$, we have that $G-S$ is locally irregular. In all cases, we note that the eventual $S$ satisfies $|S| \leq n-1$: indeed, $S$ was constructed by considering every vertex (different from $r$) exactly once, and adding to $S$ at most one incident edge on a shortest path to $r$. Thus, $\ie(G) \leq n-1$.
\end{proof}

Adapting the above strategy, we can prove that:

\begin{theorem}\label{theorem:upper-bound-cut}
For every connected graph $G$ of order $n$, with $m$ edges and maximum degree $\Delta$, we have $\ie(G) \leq \left\lfloor \frac{m}{2} \right\rfloor + n + \Delta -2$.
\end{theorem}

\begin{proof}
Let $C=(U,V)$ be a maximum cut of $G$, \textit{i.e.}, a partition of $V(G)$ maximising the number $x$ of edges from $U$ to $V$. It is well known that $x \geq \lceil m/2 \rceil$, essentially because, for every vertex $u \in U$ (or $v \in V$), we have $d_V(u) \geq d_U(u)$ ($d_U(v) \geq d_V(v)$, resp.) -- indeed, otherwise we would contradict the maximality of $C$ when swapping $u$ from $U$ to $V$ ($v$ from $V$ to $U$, resp.). We set $E(G)=B \cup R$, where $B$ contains all edges joining $U$ and $V$, and $R$ contains all other edges (each included in either $U$ or $V$).

Set $H=G-R$. Then, $H$ is bipartite. Note that $H$ is not necessarily connected. However, since $G$ is connected, there must be two (bipartite) connected components $H_1$ and $H_2$ of $H$ that are joined by an edge $e \in R$ in $G$; then, we just move $e$ from $R$ to $B$. Since adding an edge joining two disjoint bipartite graphs results in a bipartite graph, we still have that $G-R$ is bipartite. So, by repeating this process, we end up with a bipartition $B \cup R$ of $E(G)$ such that $H=G-R$ is a connected bipartite graph, and we still have $|R| \leq \lfloor m/2 \rfloor$ (since we perhaps moved edges from $R$ to $B$ in the process).

Let now $r$ be any vertex of $H$ that is not a cut-vertex of $H$ (for instance, any last vertex encountered during the run of a BFS algorithm is a candidate as $r$). Since $r$ is not a cut-vertex of $H$, by removing any set $X$ of $d(r)-1 \leq \Delta-1$ edges incident to $r$, we obtain a bipartite graph $H'$ with $\delta(H')=1$, and Theorem~\ref{theorem:bipartite-deg1-upper} tells us that $\ie(H') \leq n-1$.

If we denote by $S$ any optimal edge-irregulator of $H'$, then note that $G-R-X-S$ is locally irregular. As a result,
\begin{align*}
    \ie(G) & \leq |R \cup X \cup S| \\
    & \leq \left\lfloor \frac{m}{2} \right\rfloor + (\Delta-1) + (n-1) \\
    & = \left\lfloor \frac{m}{2} \right\rfloor + n + \Delta -2
\end{align*}

\noindent This concludes the proof.
\end{proof}

\subsection{Trees}
We now focus on trees. We first prove that Conjecture~\ref{conjecture:1/3} holds for this family, and then provide a linear-time algorithm for trees of bounded maximum degree.

\begin{figure}[t]
    \centering
    \scalebox{0.9}{
    \begin{tikzpicture}[inner sep = .5mm]
        \node[draw, black, circle, fill] (u) at (0, 0)[label = above left:$u$]{};
        \node[draw, black, circle, fill] (v) at (0, 1)[label = below right:$v$]{};
        \node[draw, black, circle, fill] (w11) at (.5, -.5)[]{};
        \node[draw, black, circle, fill] (w21) at (-.5, -.5)[]{};
        \node[draw, black, circle, fill] (w'11) at (-.5, 1.5)[]{};
        \node[draw, black, circle, fill] (w'21) at (.5, 1.5)[]{};

        \node[draw, black, circle, fill] (w12) at (1, -.5)[]{};
        \node[draw, black, circle, fill] (w13) at (1.5, -.5)[]{};
        \node[draw, black, circle] (w22) at (-1, -.5)[]{};
        \node[draw, black, circle] (w'12) at (-1, 1.5)[]{};
        \node[draw, black, circle] (w'22) at (1, 1.5)[]{};
        \node[draw, black, circle] (w23) at (-1.5, -.5)[]{};
        \node[draw, black, circle] (w'13) at (-1.5, 1.5)[]{};
        \node[draw, black, circle] (w'23) at (1.5, 1.5)[]{};
        \draw (u) -- (v);
        \draw[dashed] (u) -- (w11);
        \draw (u) -- (w21);
        \draw (v) -- (w'11);
        \draw (v) -- (w'21);
        \draw (w12) -- (w11);
        \draw (w22) -- (w21);
        \draw (w'12) -- (w'11);
        \draw (w'22) -- (w'21);
        \draw (w12) -- (w13);
        \draw (w22) -- (w23);
        \draw (w'12) -- (w'13);
        \draw (w'22) -- (w'23);
    \end{tikzpicture}
    }
    \hspace{20pt}
    \scalebox{1}{
    \begin{tikzpicture}[inner sep = .5mm]
        \node[draw, black, circle, fill] (u) at (0, 0)[label = above left:$u$]{};
        \node[draw, black, circle, fill] (v) at (0, 1)[label = below right:$v$]{};
        \node[draw, black, circle, fill] (w11) at (.5, -.5)[]{};
        \node[draw, black, circle, fill] (w21) at (-.5, -.5)[]{};
        \node[draw, black, circle, fill] (w'11) at (-.5, 1.5)[]{};
        \node[draw, black, circle, fill] (w'21) at (.5, 1.5)[]{};

        \node[draw, black, circle] (w22) at (-1, -.5)[]{};
        \node[draw, black, circle] (w'12) at (-1, 1.5)[]{};
        \node[draw, black, circle] (w'22) at (1, 1.5)[]{};
        \draw (u) -- (v);
        \draw[] (u) -- (w11);
        \draw[dashed] (u) -- (w21);
        \draw (v) -- (w'11);
        \draw (v) -- (w'21);
        \draw[dashed] (w22) -- (w21);
        \draw (w'12) -- (w'11);
        \draw (w'22) -- (w'21);
    \end{tikzpicture}
    }
    \hspace{20pt}
    \scalebox{1}{
    \begin{tikzpicture}[inner sep = .5mm]
        \node[draw, black, circle, fill] (u) at (0, 0)[label = above left:$u$]{};
        \node[draw, black, circle, fill] (v) at (0, 1)[label = below right:$v$]{};
        \node[draw, black, circle, fill] (w11) at (.5, -.5)[]{};
        \node[draw, black, circle, fill] (w21) at (-.5, -.5)[]{};
        \node[draw, black, circle, fill] (w'11) at (-.5, 1.5)[]{};
        \node[draw, black, circle, fill] (w'21) at (.5, 1.5)[]{};
        \node[draw, black, circle, fill] (w12) at (1, -.5)[]{};
        \node[draw, black, circle, fill] (w22) at (-1, -.5)[]{};
        \node[draw, black, circle, fill] (w'12) at (-1, 1.5)[]{};
        \node[draw, black, circle, fill] (w'22) at (1, 1.5)[]{};
        \draw (u) -- (v);
        \draw[] (u) -- (w11);
        \draw[dashed] (u) -- (w21);
        \draw[dashed] (w12) -- (w11);
        \draw (v) -- (w'11);
        \draw (v) -- (w'21);
        \draw[dashed] (w22) -- (w21);
        \draw (w'12) -- (w'11);
        \draw (w'22) -- (w'21);
    \end{tikzpicture}
    }
     \caption{A representation of the case where $d(u) = d(v) = 3$ and $u, v$ are adjacent in the proof of Lemma~\ref{lemma:stars}. A white vertex may or may not exist. The dashed edges form an optimal edge-irregulator. In the leftmost figure, we do not provide an edge-irregulator, but rather fall back in the case of a subdivided star.}
    \label{fig:bistars}
\end{figure}
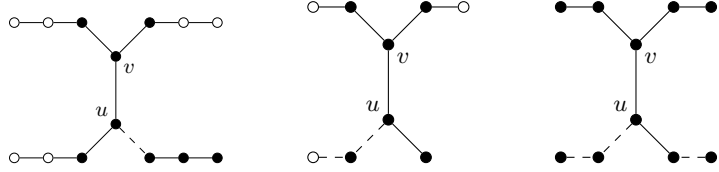

We proceed in two steps: we first deal with ``star-like'' trees, before treating the more general case. Recall that a \textit{bistar} is a tree with exactly two vertices $u$ and $v$ of degree at least~$3$ that are adjacent. See Figure~\ref{fig:bistars} for a representation of a bistar, as well as some intuition on the proof of the following lemma. For the more general case, we also need to observe a lemma concerning bridges.

\begin{lemma}\label{lemma:stars}
    If $G$ is a subdivision of a star or of a bistar, then $\ie(G)\leq\frac{1}{3}|E(G)|$. 
\end{lemma} 

\begin{proof}
Recall that the distance from an edge $e$ to a vertex $v$ is the minimum distance between the endpoints of $e$ and $v$. First assume that $G$ is a star with a center $u$ of degree at least $3$. Without loss of generality, we can assume that branches are of length at most $3$ (otherwise we remove the edge at distance 2 from the end point of such a branch to disconnect a $P_3$ an proceed on the component containing $u$). We consider $S$ to be the set of all edges at distance $2$ from $u$, which is obviously an edge-irregulator of $G$ (since all degree-two vertices of $G-S$ are adjacent to $u$ and a vertex of degree $1$). Note that $|S|$ is equal to the number of branches with exactly three edges, and all those branches are edge-disjoint, hence the result. We can now assume $G$ is a subdivided bistar other than a star, with centers $u, v$, both of degree at least $3$. 

    Assume first that $u, v$ are not adjacent in $G$ (hence, the edge between $u, v$ in the original bistar has been divided). We denote $P_{uv}$ the path between $u$ and $v$, and $u', v'$ the vertices of $P_{uv}$ adjacent to $u, v$ respectively (where $u'$ and $v'$ can be the same vertex if $d(u, v) = 2$). We find a set of edge $S_{uv}$ of $P_{uv}$ such that:
    \begin{itemize}
        \item for any two adjacent vertices $w, w'$ of $P_{uv}$ other than $u', v'$, we guarantee $w$ and $w'$ have distinct degree in $P_{uv}-S{uv}$;
        \item $uu'\notin S_{uv}$ and $vv'\notin S_{uv}$;
        \item $|S_{uv}| \leq \frac{1}{3}|E(P_{uv})|$.
    \end{itemize}
    If $|E(P_{uv})| = 0 \mod 3$ or $2 \mod 3$, then we consider $e_u, e_v$ the two edges at distance $3$ from $u, v$ respectively in $P_{uv}$. $P_{uv}-\{e_u, e_v\}$ has at most three connected components, the one containing $uu'$ (resp. $vv'$) denoted $P_{u'}$ (resp. $P_{v'}$) as well as $P$ a (possibly empty) path of length $0$ or $2 \mod 3$. We then apply the result of Theorem~\ref{theorem:paths} to obtain $S_p$ an edge-irregulator of $P$ of size at most $\frac{1}{3}|E(P)|$ and denote $S_{uv}$ the set $S_p\cup \{e_u, e_v\}$. If $|E(P_{uv})| = 1 \mod 3$, then we instead consider $e_u$ to be the edge at distance $2$ from $u$ instead, and the same process applies. It is easy to check that in both cases, all the above conditions are fulfilled.
    Finally, we consider $G_u, G_v$ each component of $G-P_{uv}$ and apply the exact same process as before on $G_u\cup P_{u'}$ and $G_v\cup P_{v'}$ (which are both subdivided stars), yielding two optimal edge-irregulators $S_u, S_v$. Note that this process does not alter the degree of neither of $u$ nor $v$. Hence $S_u\cup S_v$ is an edge-irregulator of $G$. Finally, $|S_u\cup S_v|$ it still the number of branches of length exactly $3$ or $G$. We can now assume that $u, v$ are adjacent.
    
    If $d(u)\neq d(v)$ and $d(u), d(v)\geq 3$, then we consider $G_u, G_v$ each component of $G-uv$ and apply the exact same process as before on $G_u\cup uv$ and $G_v\cup uv$ (which are both subdivided stars), yielding two optimal edge-irregulators $S_u, S_v$. Note that this process does not alter the degree of neither of $u$ nor $v$. Hence $S_u\cup S_v$ is an edge-irregulator of $G$. Finally, $|S_u\cup S_v|$ it still the number of branches of length exactly $3$ or $G$. 
    
    If $d(u) = d(v)$ and $d(u), d(v)\geq 4$, then either there is a branch of length $3$ (say attached to $u$), in which case we disconnect a $P_3$ by removing the edge connecting it to $u$, and fall back into the previous case, or all branches are of length at most $2$, and we delete an arbitrary branch. Since there are more than~$6$ edges in $v$, we are done.

    We can now assume $d(u) = d(v) = 3$. If there is a branch of length $3$, we disconnect a $P_3$ and fall back in the case of a star. If there is a branch of length $1$, denoted $uw_1$, then we denote $uw_2$ (and $w_2w_3$ if such a $w_3$ of degree $1$ exists) the edges of the other branch, and consider $S = \{uw_2, w_2w_3\}$ (if $w_2w_3$ exists). Note that if $w_2w_3$ exists, there are at least 6 edges in $G$, thus we are done. We can finally assume that all branches are of length $2$, and if we denote $w_1, w_2$ (resp. $w'_1, w'_2$) the vertices on the branches incident to $u$, then $S = \{uw_1, w_1w_2, w'_1w'_2\}$ is an edge-irregulator of $G$. Since $|E(G)| = 9$, the proof is complete.
\end{proof}

\begin{lemma}\label{lemma:bridge}
Let $S$ and $S'$ be two edge-irregulators of two graphs $G, G'$, respectively, and let $u\in V(G)$ and $v\in V(G')$. Let $H$ be the graph obtained from the disjoint union of $G$ and $G'$ with the addition of the edge $uv$. Then $I = S\cup S'\cup \{uv\}$ is an edge-irregulator of $H$.
\end{lemma}

\begin{proof}
We consider a vertex $w\in V(H)$. Assume w.l.o.g. that $w\in V(G)$. It is easy to check that the degree of $w$ in $H-I$ is the same as in $G-S$. The case $w\in V(G')$ is similar.
\end{proof}

\begin{theorem}\label{theorem:trees}
    Let $T$ be a tree. Then either $T$ is a path or $\ie(T) \leq \frac{1}{3}|E(T)|$. 
\end{theorem}

\begin{proof}
    Suppose by contradiction that $G$ is a minimum counterexample to this theorem. Then $G$ cannot be a path, and by Lemma~\ref{lemma:stars}, $G$ cannot be a subdivision of a star or of a bistar. We root $G$ at an arbitrary vertex $r$, and we consider $u$ a deepest vertex of degree $k \ge 3$ in $G$. If several exist, we choose one of the highest degree. Denote by $v$ its parent. Note that $u \neq r$, otherwise the graph would be a subdivided star. If $v = r$, we chose a neighbour of $r$ other than $u$ as the new root. 
    For $w\in V(G)$, we will denote by $T_w$ the subtree induced by $w$ and all its descendant for this rooting. 

\begin{numclaim}\label{claim no pending P3}
    Let $P=(p_1,p_2,p_3)$ be a path of $G$ with $d(p_1) = 1$ and $d(p_2) = 2$. Then $d(p_3)\geq 3$.
\end{numclaim}

\begin{proofclaim}
    Since $G$ is connected, $d(p_3) \neq 1$, otherwise we would have $T =P$, which contradicts that $G$ is not a path. Suppose by contradiction $d(p_3) = 2$, let $G' = G \setminus T_{p_3}$ and let $S'$ be an edge-irregulator of $G'$. By minimality of $G$, we have that $|S'| \le \frac{|G'|}{3}$. Let $p_3'$ be the parent of $p_3$, and consider $S = S' \cup \{p_3p_3'\}$. Then $S$ is an edge-irregulator of $G$ of size at most $\frac{|G|}{3}$, a contradiction. 
\end{proofclaim}

    It follows that, by construction, the subtree $T_u$ rooted in $u$ is a subdivided star. Let $P_i = ( v_{i,1,}, \dots, v_{i,k_i})$  for $i\in[k-1]$ be its branches. Without loss of generality, and according to Claim~\ref{claim no pending P3}, we can suppose for $1 \le i <j \le k$, $|P_i| \ge |P_j|$, and thus $2 \ge |P_1| \ge \dots \ge |P_{k}| \ge 1$.
    \begin{itemize}[noitemsep,topsep=0pt,parsep=0pt,partopsep=0pt, leftmargin=11pt]
        \item If $k \ge 4$, then by the minimality of $G$, we get that $G' = G\setminus T_u \cup \{u\}$ has an edge-irregulator $S'$ of size at most  $\frac{1}{3}|E(G')|$. If $uv \in S'$, then $S \cup S'$ is an edge-irregulator of the required size, where $S$ is an edge-irregulator of $T_u$ of size at most $\frac{|T_u|}{3}$ obtained from Lemma~\ref{lemma:stars}. Now suppose $uv \notin S'$. In $G \setminus S'$, the only remaining conflicts can be the edge $uv$. Should this happen, all the edges from $P_{k-1}$ form an edge-irregulator of $T_u$ decreasing the degree of $u$ by one, thus resolving the conflict $uv$, and $|P_i| \le \frac{|T_u|}{3}$. 
    
    \item We can now assume $k = 3$. There are thus two induced paths $P_1, P_2$ rooted at $u$. We will distinguish cases depending on the lengths of $|P_1|$ and $|P_2|$. 
    \begin{itemize}[noitemsep,topsep=0pt,parsep=0pt,partopsep=0pt, leftmargin=11pt]
        \item Assume $|P_1|=|P_2| = 1$ and let $G' = G \setminus T_u$. By minimality of $G$, the tree $G'$ is not a path and has an edge-irregulator $S'$of size $\frac{|G'|}{3}$. By construction, $S = S'\cup \{uv\}$ is an edge-irregulator of $G$. Moreover, we have $|S| = |S'| + 1$ and $|G| = |G'| + 3 $.
        \item Assume now $|P_1| = 2$ and $|P_2| \in \{1, 2\}$. Let $v'$ be the parent of $v$. If $d_G(v)=1$, then the graph is already locally irregular and we are done. If $d_G(v) = 2$ and $G'$ is a path then we are done by Lemma~\ref{lemma:stars}. Otherwise (if $d_G(v)=2$ and $G'$ is not a path),
        we consider $G' = G - T_v$ and let $S'$ be an optimal edge-irregulator of $G'$. Using Lemma~\ref{lemma:bridge}, the set $S = S' \cup \{vv'\}$ is an edge-irregulator of $G$ and has size less than $\frac{1}{3}|E(G)|$.
        
        We can now assume $d(v) \geq 3$. By the definition of $u$, any other descendant of $v$ of degree at least $3$ must be a neighbour of $v$. Moreover, by minimality of $G$, we can suppose that all its other neighbours of degree $3$ cannot fall into the previous cases. Thus, we can suppose that all its other neighbour have degree at most $3$, and if they have degree $3$, they are adjacent to a path of length two and a path of length at most $2$. Moreover, we recall that by minimality of $G$, there cannot be a path of length at least $3$ attached to $v$.
        \begin{itemize}[noitemsep,topsep=0pt,parsep=0pt,partopsep=0pt, leftmargin=11pt]
            \item Assume first $d_G(v) \geq 5$ and let $G' = G - T_v $. Consider $S'$ an optimal edge-irregulator of $G'$. Note that $G'$ can be a path, but a path has an edge-irregulator of size at most $\frac{|G'|}{3} +1$ which will be sufficient here. Let $ S = S' \cup \{vv'\}$. It is easy to check that $G-S$ is indeed locally irregular. Since exactly one edge adjacent to $v$ was removed, its resulting degree is at least $4$, all its neighbours have degree at most $3$ and the rest of $T_v$ was already locally irregular. Moreover, $T_v$ contains at least $7$ edges. Thus $|S| = |S'| +1 \le \frac{|G'|}{3} +2 \le \frac{|G|-7} {3} + 2 \le \frac{|G|}{3}$ 

            
            \item Assume $d(v) = 3$ and let $w$ be the neighbour of $v$ other than $u$ and $v'$. We consider $G' = G - T_v$. If $G'$ is a path, then $d_g(w) = 3$, as otherwise $G$ would be a bistar. In this case, consider $S'$ an edge-irregulator of $G'$ of size at most $\frac{|G'|}{3} +1$ and let $S = S' \cup \{vv'\}$. As in the previous case, since there are at least $8$ edges in $T_v$ and since $T_v$ is already locally irregular, we are done.  Otherwise, $G'$ is not a path and we can consider $S'$ an optimal edge-irregulator of $G'$ of size at most $\frac{|G'|}{3}$. If $d_G(w)\in \{1,3\}$, then $S = S' \cup \{vv'\}$ is an edge-irregulator of size at most $\frac{|G|}{3}$. Otherwise, $d_G(w) = 2$. Denote by $w'$ the child of $w$ in $T_v$. Note that $w'$ is a leaf, since we have already proved that there is no path of length at least $3$ adjacent to $v$. Thus, $S \cup \{vv', ww'\}$ is an edge-irregulator of $G$ of size at most $\frac{|G|}{3}$ as $|T_v| \ge 6$.  

            \item We can finally assume $d(v) = 4$. Let $w$ and $x$ be the two other children of $v$. W.l.o.g., suppose $d(w) \ge d(x)$. By hypothesis on $u$, we have $3 = d(u) \ge d(w) \ge d(x)$, and both $w$ and $x$ can only have paths as descendents. In case of equality in their degrees, suppose $|T_u| \ge |T_w| \ge |T_x|$.
             \begin{itemize}[noitemsep,topsep=0pt,parsep=0pt,partopsep=0pt, leftmargin=10pt]
             \item If $d(w) = 3$, let $G' = G \setminus T_v$. If $G'$ is a path, $G$ is already locally irregular. Indeed, all the vertices of degree $3$ or more are $v$ and some of its neighbours. Thus, using Claim~\ref{claim no pending P3}, and since all the neighbours or $v$ have degree at most $3$, there are no conflicts in $G$. 
                Otherwise, let $S'$ be an edge-irregulator of $G'$, and $S_x$ be an optimal edge-irregulator of $T_x$. Since $G'$ is not a path and $G$ is a minimum counter-example, we have $|S'| \le \frac{|G'|}{3}$. Also, $|S_x| \le \left \lceil \frac{|T_x|}{3} \right \rceil$ since $T_x$ is a path. Consider $S = S' \cup \{vv', vx\} \cup S_x$. By construction, $S$ is an edge-irregulator of $G$, and $|S| \le \frac{|G|}{3}$, see Figure~\ref{fig:trees}(a).
            \item If $d_G(w) = 2$ and $d_G(x) = 2$, let $G' =  G \setminus (T_u \cup T_w \cup T_x)$ and let $S'$ be an optimal edge-irregulator of $G'$. 
            If $d_{G'-S'}(v')\neq 4$, since $d_{G'-S'}(v)=4$, we have that
            $S'$ is also an edge-irregulator of $G$ and has size at most $\frac{G}{3}$. Otherwise, let $S_u$ be an optimal edge-irregulator of $T_u$. Note that $T_u$ is either a path on $4$ or on $5$ vertices, which implies $|S_u| = 1$ or $|S_u| = 2$, respectively. Consider $S = S' \cup S_u \cup \{(vu)\}$. Then $S$ is an edge-regulator of $G$, and $|S| = |S'| + 1 + |S_u| \le \frac{|G'|}{3} +1 + |S_u| \le \frac{|G|}{3}$ since $|T_v| = 8$ if $|T_u| = 4$ and $|T_v| = 9$ if $|T_u| = 5$, see Figure~\ref{fig:trees}(b).
            \item Finally, suppose $d_G(w) \le 2$, and $d_G(x) = 1$. Let $G' = G \setminus T_v$. Once again, if $G'$ is a path, $G$ is already locally irregular. Let $S'$ be an edge-irregulator of $G'$, and let $S = S' \cup \{vv', vx\}$. If $d(w) = 2$, we also add to $S$ the edge connecting $w$ to its leaf. In both cases, $S$ is an edge-irregulator of $G$ of size at most $\frac{G}{3}$, see Figure~\ref{fig:trees}(c).
            \end{itemize}


        \end{itemize}
    \end{itemize}
    \end{itemize}
    Thus, a minimal counterexample does not have a vertex of degree at least $3$, meaning it is a path, which concludes the proof.
\end{proof}

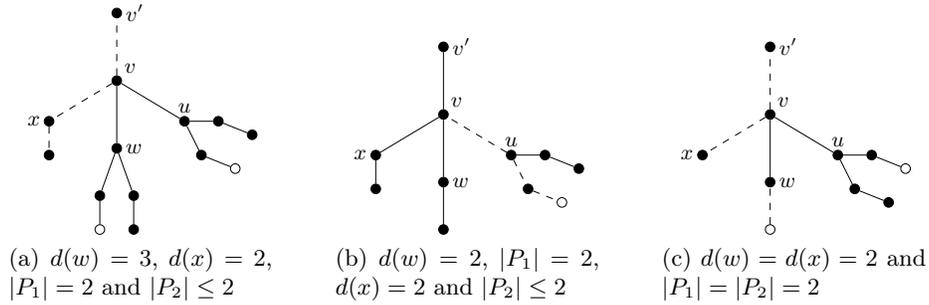
\begin{figure}[t!]
    \centering
\subfloat[$d(w) = 3$, $d(x) = 2$, $|P_1| = 2$ and $|P_2| \le 2$]{
    \scalebox{0.9}{
    \begin{tikzpicture}[inner sep = .5mm]
        \node[draw, black, circle, fill] (u) at (0, 0)[label = above:$u$]{};
        \node[draw, black, circle, fill] (w11) at (.5, 0)[]{};
        \node[draw, black, circle, fill] (w21) at (.25, -.5)[]{};
        \node[draw, black, circle, fill] (w12) at (1, -.2)[]{};
        \node[draw, black, circle] (w22) at (.75, -.7)[]{};
        \node[draw, black, circle, fill] (v') at (-1, 1.6)[label = right:$v'$]{};
        \node[draw, black, circle, fill] (v) at (-1, 0.6)[label = above right:$v$]{};
        \node[draw, black, circle, fill] (w1) at (-1, -.4)[label = right:$w$]{};
        \node[draw, black, circle, fill] (w111) at (-1.25, -1.1){};
        \node[draw, black, circle, fill] (w121) at (-.75, -1.1){};
        \node[draw, black, circle] (w112) at (-1.25, -1.6){};
        \node[draw, black, circle, fill] (w122) at (-.75, -1.6){};

        \node[draw, black, circle, fill] (w2) at (-2, 0)[label = left:$x$]{};
        \node[draw, black, circle, fill] (w211) at (-2, -.5){};

        \draw (u) -- (v);
        \draw (u) -- (w11);
        \draw (u) -- (w21);
        \draw (w11) -- (w12);
        \draw (w21) -- (w22);
        \draw[] (v) -- (w1);
        \draw[dashed] (v) -- (v');
        \draw[dashed] (v) -- (w2);
        \draw (w111) -- (w1);
        \draw (w112) -- (w111);
        \draw (w121) -- (w1);
        \draw (w122) -- (w121);
        \draw[dashed] (w2) -- (w211);
    \end{tikzpicture}
    }
    }\hfill
    \subfloat[$d(w) = 2$, $|P_1| = 2$, $d(x) = 2$ and $|P_2| \le 2$]{
    \scalebox{0.9}{
    \begin{tikzpicture}[inner sep = .5mm]
        \node[draw, black, circle, fill] (u) at (0, 0)[label = above:$u$]{};
        \node[draw, black, circle, fill] (w11) at (.5, 0)[]{};
        \node[draw, black, circle, fill] (w21) at (.25, -.5)[]{};
        \node[draw, black, circle, fill] (w12) at (1, -.2)[]{};
        \node[draw, black, circle] (w22) at (.75, -.7)[]{};
        \node[draw, black, circle, fill] (v') at (-1, 1.6)[label = right:$v'$]{};
        \node[draw, black, circle, fill] (v) at (-1, 0.6)[label = above right:$v$]{};
        \node[draw, black, circle, fill] (w1) at (-1, -.4)[label = right:$w$]{};
        \node[draw, black, circle, fill] (w111) at (-1, -1.1){};
        \node[draw, black, circle, fill] (w2) at (-2, 0)[label = left:$x$]{};
        \node[draw, black, circle, fill] (w211) at (-2, -.5 ){};

        \draw[dashed] (u) -- (v);
        \draw (u) -- (w11);
        \draw[dashed] (u) -- (w21);
        \draw (w11) -- (w12);
        \draw[dashed] (w21) -- (w22);
        \draw (v) -- (w1);
        \draw (v) -- (v');
        \draw (v) -- (w2);
        \draw (w111) -- (w1);
        \draw (w2) -- (w211);
    \end{tikzpicture}
    }
    }\hfill    
    \subfloat[$d(w) = d(x) = 2$ and $|P_1| = |P_2| = 2$]{
    \scalebox{0.9}{
    \begin{tikzpicture}[inner sep = .5mm]
        \node[draw, black, circle, fill] (u) at (0, 0)[label = above:$u$]{};
        \node[draw, black, circle, fill] (w11) at (.5, 0)[]{};
        \node[draw, black, circle, fill] (w21) at (.25, -.5)[]{};
        \node[draw, black, circle] (w12) at (1, -.2)[]{};
        \node[draw, black, circle, fill] (w22) at (.75, -.7)[]{};
        \node[draw, black, circle, fill] (v') at (-1, 1.6)[label = right:$v'$]{};
        \node[draw, black, circle, fill] (v) at (-1, 0.6)[label = above right:$v$]{};
        \node[draw, black, circle, fill] (w1) at (-1, -.4)[label = right:$w$]{};
        \node[draw, black, circle] (w111) at (-1, -1.1){};
        \node[draw, black, circle, fill] (w2) at (-2, 0)[label = left:$x$]{};

        \draw (u) -- (v);
        \draw (u) -- (w11);
        \draw[] (u) -- (w21);
        \draw (w11) -- (w12);
        \draw (w21) -- (w22);
        \draw[] (v) -- (w1);
        \draw[dashed] (v) -- (v');
        \draw[dashed] (v) -- (w2);
        \draw[dashed] (w111) -- (w1);
    \end{tikzpicture}
    }
    }
    \caption{A representation of some cases of the proof of Theorem~\ref{theorem:trees}. The white vertices may or may not exist. The dashed edges are the edges in $S$.}\label{fig:trees}
\end{figure}

Theorem~\ref{theorem:trees} proves that Conjecture~\ref{conjecture:1/3} holds for trees. For this family we can also compute the exact value of $\ie(G)$ with an efficient algorithm. 

\begin{theorem}\label{theorem:algo-trees}
If $G$ is a tree of order $n$ and maximum degree $\Delta$, then computing $\ie(G)$ can be done in time $O(n\Delta^3)$.
\end{theorem}

\begin{proof}
    Intuitively, we will compute recursively, for each vertex $u$, two arrays of size at most $\Delta$ containing the minimum number of edges in an edge-irregulator of the subtree induced by $u$, with or without the edge to its parent, for any possible degree for $u$. We fix a root for $G$, and denote $T_u$ the subtree induced by $u$, and $p$ its parent. We denote $U^u$ the array that will contain, at the $i$th position, the minimum number of edges in an edge-irregulator $S$ of $T_u$ under the assumption that $u$ has degree $i$ in $T_u\cup (pu)$. Similarly, we denote $D^u$ the array that will contain, at the $i$th position, the minimum number of edges in an edge-irregulator $S$ of $T_u$ under the assumption that $u$ has degree $i$ in $T_u$. Note that if $u$ is a leaf, $U^u$ and $D^u$ are both the array of size one containing a zero. 

    Now, we suppose $u$ is not a leaf and denote $v_1, \dots, v_k$ its children. We first compute the values of $D^u$. For any $1\leq \delta \leq k$ which will be the prescribed degree of $u$ (and the value of $D^u(\delta)$), for any $1\leq i \leq k$, we denote $in(i)$ the smallest value of $U^{v_i}$ other than $U^{v_i}(\delta)$. Similarly, we denote $out(i)$ the smallest value of $D^{v_i}$. We claim that $D^u(\delta)$ is equal to $k-\delta + \Sigma_{i\in S_1} in(i) + \Sigma_{i\in S_2} out(i)$ where $S_1$ is the subset of $\{1, \dots, k\}$ containing the values of the $\delta$ indices $i$ such that $in(i)-out(i)$ is minimal, and $S_2$ contains the other indices. We compute each value of $Up^u$ in the exact same fashion (if $u$ is not the root), except that $S_1$ will contain $\delta-1$ indexes instead of $\delta$. Note that if we denote $r$ the root of $G$, then $\ie(G)$ is the smallest value of $D^r$.

\noindent\textbf{Running time.} We now prove that this algorithm runs in time $O(n\Delta^3)$. For each vertex, and for each degree prescription $\delta$, computing the $in$ and $out$ functions, as well as the values $in-out$, can be done in time $O(\Delta)$. Searching the $\delta$ smallest values of $in-out$ can be done in time $O(\Delta^2)$, hence the final complexity.

\noindent\textbf{Correctness.} Note that reaching degree $\delta$ for $u$ implies choosing $\delta$ children such that the edge between $u$ and said child in not deleted. In particular, if we denote by $S_1$ the set of children adjacent to $u$ and $S_2$ the set of children not adjacent to $u$, then the size of an optimal edge-irregulator of $T_u$ is smaller than $\Sigma_{i\in S_1} in(i) + \Sigma_{i\in S_2} out(i) + d_G(u) - \delta= d_G(u) - \delta + \Sigma_{i\in S_1} in(i) + \Sigma_{i\in S_2} out(i) = d_G(u) - \delta + \Sigma_{i} in(i) - \Sigma_{i\in S_2} in(i) + \Sigma_{i} out(i) - \Sigma_{i\in S_1} out(i) = d_G(u) - \delta + \Sigma_{i\in S_1} (in(i) - out(i)) + \Sigma_{i} out(i)$.

Since the values $\Sigma_{i} out(i)$, $d_G(u)$ and $\delta$ are fixed, choosing the minimal value for $\Sigma_{i\in S_1} (in(i) - out(i))$ yields an optimal edge-irregulator of $T_u$.
\end{proof}

\subsection {Complete graphs}

For any $n \geq 3$, we denote by $K_n$ the complete graph of order $n$. Complete graphs are an interested case for our problem, as computing $\ie(K_n)$ is the same as finding the size of locally irregular graphs of order $n$ with the most edges.

In what follows, we are solving the question for particular values of $n$, while we only provide partial results for the remaining ones. Precisely, we focus on the case where $n$ is a \textit{triangular number}, \textit{i.e.}, of the form $\frac{k(k+1)}{2}$ for some $k$. For any $k \geq 1$, we denote by $t_k$ the $k^{th}$ triangular number, being the sum of all integers from $1$ up to $k$.

\begin{lemma}\label{lemma:triangular-li}
If $n=t_k$ is a triangular number, 
then the maximum number of edges in a locally irregular graph of order $n$ is exactly $m_k:=\frac{k(k+1)(k-1)(3k+2)}{24}$.
\end{lemma}

\begin{proof}
Let $G$ be a locally irregular graph of order $n=t_k$.
First, note that, for every $i \in \{1,\dots,k\}$, there are at most $i$ vertices of degree $n-k$ in $G$, otherwise, by the pigeonhole principle, two of them would be adjacent. So, any locally irregular graph of order $n$ having the most edges has one vertex of degree $n-1$, two vertices of degree $n-2$, and so on.
Because $n$ is a triangular number, such a graph $T_k$ can be constructed for all $k \geq 1$; we prove this by induction on $k$.

The graph $T_1$ is the graph containing a single vertex. Then, for any $i \geq 1$ such that $T_i$ has been constructed, we obtain $T_{i+1}$  from $T_i$ by adding $i+1$ independent vertices adjacent to all the vertices of $T_i$. It is easy to see that for any $k \geq 1$, the graph $T_k$ indeed has the claimed degree sequence. To fully prove the claim, it remains to express the number of edges of any $T_k$. The starting point is the handshaking lemma, which states that in any graph the number of edges is twice the sum of the degrees:
\begin{align*}
    2|E(T_k)| &= \underset{i=1}{\overset{k}{\sum}} i(n-i)\\
    &= n\underset{i=1}{\overset{k}{\sum}} i - \underset{i=1}{\overset{k}{\sum}} i^2\\
    &= \left (\frac{k(k+1)}{2}\right)^2 - \frac{k(k+1)(2k+1)}{6}  \tag*{\mbox{\textit{since $n = \frac{k(k+1)}{2}$}}}\\
    &= \frac{k(k+1)}{2} \left ( \frac{k(k+1)}{2} - \frac{2k+1}{3} \right ) \\
    &= \frac{k(k+1)}{2} \frac{3k^2-k-2}{6}\\
    &= \frac{k(k+1)(k-1)(3k+2)}{12}
\end{align*}
This concludes the proof.
\end{proof}

\begin{theorem}\label{theorem:complete-triangular}
For every $k \geq 1$ and $n = t_k$, we have $\ie(K_n)=|E(K_n)|-m_k$.
\end{theorem}

Even though Theorem~\ref{theorem:complete-triangular} is tight, it cannot be generalised easily to complete graphs whose order is not a triangular number. However, since $t_{k+1} - t_k = k = \Theta(\sqrt{t_k})$ for every $k \geq 1$, Lemma~\ref{lemma:triangular-li} yields both a lower and an upper bound on the size of locally irregular graphs of order $n$ having the most edges, when $n$ is not triangular: namely, they have $\frac{n^2}{2} - O(n^{3/2})$ edges. And from this, we proved that Conjecture~\ref{conjecture:1/3} holds for all complete graphs. 

\section{Conclusion} \label{sec:ccl}

In this work, we extended the study of edge-irregulators from~\cite{FMT24}. Apart from focusing on the algorithmic aspects of this problem, we also derived general bounds and analysed behaviours of optimal edge-irregulators, leading to Conjecture~\ref{conjecture:1/3}, which we verified for several graph classes.

Determining $\ie(G)$ for a graph $G$ remains challenging due to the intricate properties and elusive structure of maximally edge-dense locally irregular graphs. Let us denote by $\mathcal{L}_n$ (for some $n \geq 3$) the class of all locally irregular graphs of order $n$ with the most edges possible. In light of Lemma~\ref{lemma:triangular-li}, a key open question concerns the size and structural properties of $\mathcal{L}_n$ in the case where $n$ is a non-triangular number.

We finally observe that there is room for improvement for many of our results. This could lead to further work on the topic. In particular, we find the following questions appealing.

\begin{itemize}
    \item Refining the static parameter $\conf$ (Section~\ref{sec:structure}) to dynamically capture evolving edge conflicts.

    \item Determining $\ie$ for complete graphs and tightening bounds for cubic graphs.

    \item Extending Theorem~\ref{theorem:subdiv-1/3} to subcubic graphs. Also, for $k$-degenerate graphs $G$ of order $n$, we could get an upper bound on $\ie(G)$ of about $(k-1)n+\left\lfloor \frac{n}{3} \right\rfloor$ through exploiting similar ideas to the ones used in this paper. 
    \end{itemize}

It is also interesting to further study the algorithmic complexity of determining $\ie(G)$ for a given graph $G$. In particular, the authors of~\cite{FMT22} proved that the vertex counterpart of the problem is in \FPT~when parameterized by the treewidth and the maximum degree. Their algorithm should adapt in an obvious way to our context. Thus, we prefer to not include such a result here.
    
We finish by commenting on the \FPT~algorithm we present here w.r.t. the vertex cover number. First, observe that there is still room for improvement in its running time. Indeed, we are only aware of a $2^{o(vc)}n^{O(1)}$ lower bound (based on the ETH) for this problem (stemming as a direct corollary of the $\W[1]$-hardness w.r.t. the solution size proof presented in~\cite{FMT24}). Lastly, we believe that this algorithm could be generalised to work w.r.t. the neighbourhood diversity of the input graph. However, in order to achieve this, it would be mandatory to completely determine $\ie$ for complete graphs.


\bibliographystyle{splncs04}

\bibliography{main}

\end{document}